# Non-crossing Knight's Tour in 3-Dimension

Awani Kumar. 2/11-C, Vijayant Khand, Gomti Nagar, Lucknow 226010. INDIA

E-mail: awanieva@gmail.com

**Abstract**

Non-crossing knight's tour in 3-dimension is a new field of research. The author has shown its possibility in small cuboids and in cubes up to 8x8x8 size. It can also be extended to larger size cubes and cuboids. The author has achieved jumps of length 15, 46, 88, 159, 258 and 395 in cubes of size 3x3x3, 4x4x4, 5x5x5, 6x6x6, 7x7x7 and 8x8x8 respectively. This amounts to covering 59%, 73%, 71%, 74%, 76% and 77% cells in these cubes.

## Introduction

The problem of knight's tour on a square board is almost as old as the game itself but the non-crossing (or the non-intersecting) tour problem is a recent one. Here, the knight needs not to visit all the squares but care should be taken to find the longest path on a given board without visiting a square twice or crossing its own path. Here 'path' consists of straight lines drawn between the centre of the starting and ending square of every jump. Fig.1 and Fig.2 are examples of non-crossing tours of length 5 on a 4x4 board. The former is an open tour and the later, a reentrant tour. Lines show the non-crossing path. Basically, knight is a three dimensional piece and its move (0, 1, 2) has an aesthetic appeal being the first three whole number. However, traditional study of knight's tour has been mostly confined to 2-dimension. Perusal of literature reveals that the problem of non-crossing path on the 8×8 board was solved in 1930, with an open path of 35 moves by T. R. Dawson and with a closed path of 32 moves by the Romanian chess problemist Wolfgang Pauly [1]. These results were reported without a diagram of Pauly's result. Later, Murray [2] showed the diagram in his unpublished manuscript. The knight problem for small rectangular boards was rediscovered by Yarbrough [3]. Some of his results were improved on in letters in the same journal 1969 (vol.2, nr.3, pp.154-157) by R. E. Ruemmler (7×8 and 5×9 to 9×9), D. E. Knuth (5×6, 6×6, 7×8, 8×8, confirming the Dawson/Pauly results, and 5×9) and M. Matsuda (6×6, 6×8, 5×9, 7×9 and 9×9). Jelliss [4] and Merson [5] have looked into non-crossing tours by fairy pieces and on larger boards.

More recently, Awani Kumar [6] looked into the problem of knight's tour in 3-dimension but it was mostly confined to tours having magic properties. Now, the author proposes to look into the problem of non-crossing knight's tour in 3-dimension. Here 'path' consists of straight lines drawn between the centre of the starting and ending cells of every jump.

**Non-crossing tours in cuboids**: In 2-dimension, the readers can easily see that 2x3 is the smallest rectangle in which a knight can move. In 3-dimension, 2x2x3 is the smallest cuboid in which a non-crossing knight tour of length 4 is possible. The readers can visualize it in 3-dimension by stacking the 2x3 rectangles, one over the other, in alphabetical order as shown in Fig.3. Fig.4 shows a non-crossing knight's tour of length 8 in a 2x3x3 cuboid. It is interesting to note that Fig.3 and Fig.4 are non-crossing closed knight’s tour, that is, their starting and ending cells are the same. Such tours are rare. Fig.5, Fig.6 and Fig.7 are tours in cuboids of size 2x4x4, 3x3x4 and 3x4x4 and having length 14, 20 and 27 respectively.

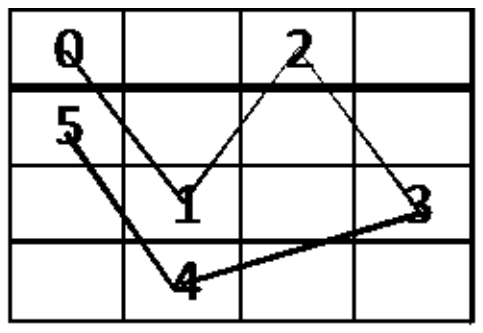

Fig.1. Open tour

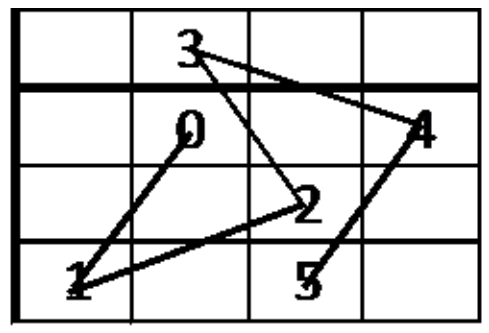

Fig.2. Reentrant tour

A

| | | |
|---|---|---|
| 0,4 | | |
| | | 1 |

B

| | | |
|---|---|---|
| | | 3 |
| 2 | | |

Fig.3. Non-crossing closed knight's tour in 2x2x3 cuboid

A

| | | |
|---|---|---|
| | 1 | 6 |
| | | |
| 0,8 | | 2 |

B

| | | |
|---|---|---|
| 7 | 4 | |
| | | |
| 3 | | 5 |

Fig.4. Non-crossing closed knight's tour in 2x3x3 cuboid

A

| | | | |
|---|---|---|---|
| 0 | | 2 | |
| 5 | 12 | | |
| | 1 | | 3 |
| | 4 | | 14 |

B

| | | | |
|---|---|---|---|
| | | 9 | |
| 8 | | | 13 |
| | | 7 | 10 |
| 6 | 11 | | |

Fig.5. Non-crossing knight's tour in 2x4x4 cuboid

A

| | | | |
|---|---|---|---|
| 1 | 4 | 17 | |
| | | 0 | 3 |
| | 2 | 19 | 14 |

B

| | | | |
|---|---|---|---|
| | 7 | | 15 |
| | | 9 | 6 |
| 8 | 5 | | |

C

| | | | |
|---|---|---|---|
| | 16 | 11 | |
| 10 | | 18 | 13 |
| | 12 | | 20 |

Fig.6. Non-crossing knight's tour in 3x3x4 cuboid

A

| | | | |
|---|---|---|---|
| 3 | 6 | | 26 |
| 20 | 9 | | 7 |
| | | 13 | 26 |
| | | 18 | |

B

| | | | |
|---|---|---|---|
| 22 | | 2 | 5 |
| | 0 | 23 | |
| 12 | 15 | | 1 |
| 19 | | | 16 |

C

| | | | |
|---|---|---|---|
| | 4 | | |
| 10 | | 8 | 25 |
| 21 | | 11 | 14 |
| | 17 | 24 | 27 |

Fig.7. Non-crossing knight's tour in 3x4x4 cuboid

**Non-crossing tours in cubes**: In 2-dimension, 3x3 is the smallest square in which a knight can move. In 3-dimension, 3x3x3 is the smallest cube in which a non-crossing knight tour of length 15 is possible as shown in Fig.8. Fig.9, Fig.10, Fig.11, Fig.12, and Fig.13 are tours in cubes of size 4x4x4, 5x5x5, 6x6x6, 7x7x7 and 8x8x8 having length 46, 88, 159, 258 and 395 respectively. Readers are encouraged to look for longer tours in these cubes and cuboids.

| | | |
|---|---|---|
| 9 | | 3 |
| | | 8 |
| 1 | 10 | |

A

| | | |
|---|---|---|
| 2 | 13 | |
| 7 | | |
| 12 | | 14 |

B

| | | |
|---|---|---|
| | 4 | 15 |
| 0 | | 6 |
| 5 | | 11 |

C

Fig.8. Non-crossing reentrant knight's tour in 3x3x3 cube

| | | | |
|---|---|---|---|
| 3 | | 5 | 40 |
| 10 | 41 | 30 | |
| 29 | 20 | 39 | 32 |
| | 31 | 24 | |

A

| | | | |
|---|---|---|---|
| 8 | | | 33 |
| | 22 | 9 | 12 |
| 2 | | 18 | 21 |
| 23 | 42 | 37 | 26 |

B

| | | | |
|---|---|---|---|
| 1 | 4 | 15 | 6 |
| 14 | 11 | 38 | |
| 19 | 16 | | 34 |
| 28 | 25 | | |

C

| | | | |
|---|---|---|---|
| | 7 | 45 | 35 |
| 45 | | 13 | |
| 0 | 43 | 36 | 17 |
| | 46 | 27 | |

D

Fig.9. Non-crossing knight's tour in 4x4x4 cube

| | | | | |
|---|---|---|---|---|
| 9 | | 3 | 70 | |
| 60 | 69 | 12 | | 16 |
| 27 | | 31 | 68 | 71 |
| | 67 | 72 | 37 | |
| 51 | | | | 49 |

A

| | | | | |
|---|---|---|---|---|
| | 5 | 10 | | 2 |
| 11 | 66 | 15 | 86 | |
| | 85 | 28 | | 32 |
| | 36 | 39 | 84 | 87 |
| 58 | | 50 | 53 | |

B

| | | | | |
|---|---|---|---|---|
| 61 | 8 | 77 | 4 | |
| 20 | 13 | 46 | 17 | 78 |
| 59 | 26 | 75 | 30 | 47 |
| | 45 | | 41 | 38 |
| | 52 | 73 | 48 | 55 |

C

| | | | | |
|---|---|---|---|---|
| | 65 | 6 | 1 | |
| | 34 | 19 | 14 | |
| 62 | 29 | 24 | 33 | 80 |
| 35 | 40 | 43 | | |
| | 57 | 54 | 83 | 88 |

D

| | | | | |
|---|---|---|---|---|
| 7 | 0 | | 22 | 79 |
| | 21 | 76 | | 18 |
| 25 | 64 | | 82 | 23 |
| 44 | | 74 | | 42 |
| 63 | | | 56 | 81 |

E

Fig.10. Non-crossing knight's tour in 5x5x5 cube

| | | | | | |
|---|---|---|---|---|---|
| 141 | | 1 | | 13 | 52 |
| 32 | | 142 | 53 | 30 | |
| 43 | 140 | 47 | 158 | 51 | 136 |
| 152 | 57 | | 143 | 54 | 113 |
| 139 | 78 | | | 137 | 80 |
| | 97 | 138 | | | 99 |

A

| | | | | | |
|---|---|---|---|---|---|
| 0 | | 12 | 157 | 16 | 135 |
| 153 | 156 | 31 | 28 | | 112 |
| 150 | | 44 | 155 | 48 | |
| 59 | 154 | 65 | 56 | | |
| 76 | | 82 | 79 | 130 | |
| 95 | | 101 | 98 | | 114 |

B

| | | | | | |
|---|---|---|---|---|---|
| 149 | 2 | 17 | 10 | | 14 |
| | 33 | 26 | 159 | 18 | 29 |
| | 42 | | 96 | 129 | 50 |
| | 63 | 58 | 71 | 66 | 55 |
| 151 | 72 | 77 | 144 | 81 | 84 |
| | 103 | 96 | 107 | 100 | |

C

| | | | | | |
|---|---|---|---|---|---|
| | 11 | 4 | 15 | 8 | |
| | 24 | | 20 | 27 | 134 |
| 148 | 45 | 40 | 49 | 36 | 111 |
| | 60 | 69 | 64 | 133 | |
| | 75 | 86 | 83 | 110 | 115 |
| | 94 | 109 | 102 | 131 | 106 |

D

| | | | | | |
|---|---|---|---|---|---|
| 3 | | 9 | 6 | | |
| 34 | | 22 | 25 | 128 | 19 |
| 41 | 146 | 35 | 38 | 89 | 116 |
| | 91 | 62 | 145 | 70 | 67 |
| 147 | | 73 | 90 | 85 | 88 |
| 92 | | 104 | | 108 | |

E

| | | | | | |
|---|---|---|---|---|---|
| | 5 | | 123 | | 7 |
| 23 | 122 | | | 21 | 124 |
| 120 | 39 | | 125 | | 37 |
| 61 | | 121 | 68 | 127 | |
| 74 | 119 | 126 | 87 | 132 | 117 |
| | | 93 | 118 | 105 | |

F

Fig.11. Non-crossing knight's tour in 6x6x6 cube

| | | | | | | |
|---|---|---|---|---|---|---|
| | 240 | | 26 | 251 | | 15 |
| | | 58 | 241 | 54 | 249 | 50 |
| 239 | 242 | 83 | 250 | 61 | | 65 |
| | | 104 | | 86 | | 248 |
| 243 | 220 | 217 | 118 | 247 | 202 | 129 |
| 150 | | 146 | 245 | 142 | | 190 |
| 171 | 244 | 219 | | 169 | 246 | 165 |

A

| | | | | | | |
|---|---|---|---|---|---|---|
| 238 | 27 | 24 | | 14 | | 198 |
| 59 | | 55 | | 51 | | |
| 84 | | 60 | 81 | 64 | 199 | |
| 103 | | 85 | 106 | 253 | 88 | 189 |
| | 117 | 120 | | 130 | | 200 |
| 233 | 11 | 149 | | 145 | | 141 |
| | 221 | 170 | 173 | 166 | 201 | 192 |

B

| | | | | | | |
|---|---|---|---|---|---|---|
| 29 | | 35 | 12 | 25 | 16 | |
| 36 | 57 | | 53 | 252 | 49 | 188 |
| 79 | 82 | | 62 | 187 | 66 | |
| 236 | 105 | 216 | 87 | 108 | 203 | 90 |
| 115 | 222 | 109 | 132 | 119 | 128 | 191 |
| | 151 | 218 | 147 | 254 | 143 | |
| 175 | 172 | | 168 | | 164 | |

C

| | | | | | | |
|---|---|---|---|---|---|---|
| | 33 | 28 | 23 | 18 | 13 | |
| 237 | 40 | 215 | 56 | 47 | 52 | 197 |
| 8 | 75 | 80 | | 68 | 63 | |
| 11 | 102 | 107 | 98 | 89 | 94 | |
| 234 | 111 | 116 | 121 | 126 | 131 | |
| 5 | 148 | 153 | 144 | 135 | 140 | 193 |
| 232 | 179 | 174 | | 162 | 167 | |

D

| | | | | | | |
|---|---|---|---|---|---|---|
| 9 | 30 | 21 | 34 | 11 | 208 | 17 |
| | 37 | 10 | 41 | 186 | 45 | 48 |
| 235 | 78 | | 74 | | 70 | 67 |
| 6 | 99 | 214 | 95 | 184 | 91 | 194 |
| 229 | 114 | 123 | 110 | 133 | 204 | 127 |
| 152 | 223 | 134 | 155 | | 137 | 256 |
| | 176 | | 180 | 255 | 160 | 163 |

E

| | | | | | | |
|---|---|---|---|---|---|---|
| | 227 | 32 | | 22 | 19 | |
| 7 | | 39 | 46 | 43 | | |
| 228 | 3 | 76 | 69 | 72 | 207 | 196 |
| 101 | | 97 | 2 | 93 | | |
| 4 | 225 | 112 | 211 | 122 | 125 | 258 |
| 231 | 154 | 1 | 136 | 157 | | 139 |
| 0 | | 178 | 161 | 182 | 205 | 158 |

F

| | | | | | | |
|---|---|---|---|---|---|---|
| 31 | | | 20 | | | |
| 38 | | 42 | | | 209 | 44 |
| 77 | 226 | 73 | 210 | 185 | | 71 |
| 230 | | 100 | | 96 | | 92 |
| 113 | | 213 | 124 | 183 | 206 | 195 |
| | | 156 | | 138 | | |
| 177 | 224 | 181 | 212 | 159 | | 257 |

G

Fig.12. Non-crossing knight's tour in 7x7x7 cube

| | | | | | | | |
|---|---|---|---|---|---|---|---|
| 12 | 319 | | 317 | 14 | 315 | | 17 |
| 53 | | | 50 | | 312 | | 48 |
| 320 | 95 | 318 | 81 | 316 | | 314 | 83 |
| 119 | 390 | | 134 | 313 | | 331 | 136 |
| 196 | 321 | 326 | 193 | 330 | | | 191 |
| 323 | | 165 | | 327 | 162 | | 150 |
| 250 | | 322 | 325 | 248 | 329 | | 227 |
| 285 | 324 | | 328 | 283 | | | 280 |

A

| | | | | | | | |
|---|---|---|---|---|---|---|---|
| 340 | 13 | 10 | | 16 | | 5 | |
| | 43 | 52 | 339 | 36 | 49 | 332 | 303 |
| 341 | 338 | 97 | 94 | 333 | 82 | 79 | 4 |
| 344 | 335 | 118 | | 132 | 135 | 148 | 305 |
| 337 | 186 | 195 | 334 | 189 | 192 | 3 | |
| 166 | 391 | 336 | 163 | 170 | 149 | 160 | 293 |
| | 254 | 249 | 246 | | 228 | | 0 |
| 346 | | 284 | 287 | | 281 | 292 | |

B

| | | | | | | | |
|---|---|---|---|---|---|---|---|
| 342 | 11 | 26 | 15 | 8 | 311 | 18 | |
| | 54 | | 44 | 51 | 62 | 47 | 310 |
| 96 | 389 | 92 | | 80 | 309 | 84 | |
| | 120 | 133 | 116 | 147 | | 137 | 308 |
| 348 | 197 | 202 | 187 | 194 | 205 | 190 | |
| 345 | 164 | 169 | | 161 | 172 | 151 | 306 |
| | 251 | | 255 | 230 | 247 | 226 | |
| | 286 | 271 | 282 | 289 | | 279 | |

C

| | | | | | | | |
|---|---|---|---|---|---|---|---|
| 27 | 360 | 9 | 24 | | 6 | 21 | 302 |
| 358 | 57 | 42 | 361 | 60 | 45 | | 63 |
| 343 | 98 | 359 | 90 | 93 | 78 | 87 | 304 |
| | 117 | 122 | 131 | 362 | 145 | | 139 |
| 353 | 200 | 185 | 218 | 203 | 188 | 363 | 206 |
| 180 | 167 | | 171 | 184 | 159 | 2 | 153 |
| 347 | 392 | 253 | 240 | 245 | 224 | 229 | 294 |
| 270 | | 288 | 273 | | 291 | 276 | 1 |

D

| | | | | | | | |
|---|---|---|---|---|---|---|---|
| | 25 | 40 | 7 | 22 | 37 | | 19 |
| 55 | | 73 | 58 | 41 | 76 | 61 | |
| 354 | 91 | 100 | 77 | 88 | 103 | | 85 |
| 121 | 388 | 115 | 146 | 129 | 138 | 143 | 300 |
| 198 | 393 | 216 | 201 | | 219 | 204 | 307 |
| 349 | 386 | 181 | 168 | 173 | 152 | 157 | 296 |
| 252 | 239 | 394 | 243 | 256 | 233 | 364 | 225 |
| 395 | 272 | 257 | 290 | 275 | 260 | | 278 |

E

| | | | | | | | |
|---|---|---|---|---|---|---|---|
| 357 | 28 | 23 | 38 | 31 | 20 | 35 | |
| | 71 | 56 | 67 | 74 | 59 | 64 | 301 |
| 99 | 112 | 89 | 102 | 109 | 86 | 105 | |
| | 123 | 130 | 113 | 144 | 127 | 140 | |
| 355 | 214 | 199 | 210 | 217 | | 207 | 298 |
| 352 | 179 | | 183 | 158 | 175 | 154 | |
| 237 | | 241 | | 233 | 244 | 223 | |
| 350 | 269 | 274 | 259 | 266 | 277 | 262 | 295 |

F

| | | | | | | | |
|---|---|---|---|---|---|---|---|
| | 39 | 30 | | 36 | 33 | | |
| | | 69 | 72 | | 66 | 75 | |
| 356 | 101 | 110 | | 104 | 107 | | 299 |
| | 114 | 125 | 128 | | 142 | 367 | |
| | 387 | 212 | 215 | | 209 | 220 | |
| | 182 | 177 | 174 | | 156 | | |
| | 385 | 238 | 235 | 242 | 221 | 232 | 297 |
| | 258 | 267 | 384 | 261 | 264 | 365 | |

G

| | | | | | | | |
|---|---|---|---|---|---|---|---|
| 29 | 372 | | 32 | | 370 | | 34 |
| 70 | | | 371 | 68 | | 368 | 65 |
| 111 | | 373 | 108 | 369 | | | 106 |
| 124 | 375 | | | 126 | | | 141 |
| 213 | | | 374 | 211 | 382 | | 208 |
| 178 | | 376 | 383 | 176 | 379 | 366 | 155 |
| 351 | 236 | | 378 | | 234 | 381 | 222 |
| 268 | 377 | | 265 | 380 | | | 263 |

H

Fig.13. Non-crossing knight's tour in 8x8x8 cube

**Conclusion**: About non-crossing tours in 2-dimension, Meyrignac [7] states, “Strangely, this problem did not receive much consideration, although it's **much harder** to solve than the normal knight tour problem! For example: the knight tour problem can be solved via the Warnsdorf heuristic or a recent divide-and-conquer algorithm among other methods, and we still don't know any heuristic on the uncrossing variant.” In 3-dimension, the problem is even **more hard**. The author has shown the possibility of non-crossing knight’s tour in 3-dimension and has set the ball rolling (or let loose the knight!) by venturing into small cubes and cuboids. The percentage of cells covered has been increasing with the increase in the size of the cube and the author has achieved 77% coverage in an 8x8x8 cube. What can be the maximum percentage in a cube? What can be the maximum length in mxnxk cuboid? What about non-crossing knight’s tour in 4-dimension? We count on readers to pursue this search. From little acorns, let us grow mighty oaks.

**Acknowledgement:** The author is grateful to G. P. Jelliss for providing photocopy of Ref.2, a priceless gem, and for his wonderful site on knight's tour at http://www.ktn.freeuk.com. The author is also grateful to Kirk Bresniker for inquiring about the possibility of non-crossing knight's tour in 3-dimension.